\documentclass[11pt,a4paper]{article}
\usepackage[latin1]{inputenc}
\usepackage[english]{babel}
\usepackage{amsmath}
\usepackage{amsfonts}
\usepackage{amssymb}
\usepackage{graphicx}
\usepackage{url}
\author{Moritz L. Sümmermann\thanks{moritz@suemmermann.de, Institute for Mathematics Education, University of Cologne, Germany}}
\title{Knotted Portals in Virtual Reality}
\usepackage[autostyle=true,english=american]{csquotes}
\usepackage[left=2cm,right=2cm,top=2cm,bottom=2cm]{geometry}
\usepackage{microtype}
\usepackage{xspace}
\usepackage{float}
\usepackage{hyperref}
\usepackage{pgfplots}

\usepackage{subcaption}

\begin{document}
	\maketitle
	\begin{abstract}
		\textsc{KnotPortal} is a software for the visualization of branched covers of knots based on an idea by Bill Thurston~\cite{Thurston2012}. It imagines knots made of a magical material which \enquote{rips the universe apart}, leading to the creation of portals to other worlds. This makes possible the visualization of three-manifolds constructed through gluing of different sheets along the knot as a branching curve. To recreate the experience of \enquote{stepping through the knot} described by Thurston, our implementation allows users to explore these knotted portals in virtual reality using a head-mounted device with room-tracking. Users not in possession of such a device can alternatively use the software on a normal computer screen and with keyboard controls.
		
		This article gives a short introduction into branched coverings and the history of branched covers of knots as well as the mathematical background to the ideas described by Thurston and used in the software. It also provides examples of branched coverings and the associated deck transformation groups, which are required as input for \textsc{KnotPortal}.
		
		\textsc{KnotPortal} can be used to enable students to learn about knots, gluing, (branched) covers, or just to have a fun looking at portals and knots. It is open-source and available for free download at the website of the imaginary foundation at \url{https://imaginary.org/program/knotportal}.
	\end{abstract}
	
	\section{Introduction}
\begin{figure}[H]
			\minipage[t]{0.48\textwidth}
			\includegraphics[width=\linewidth]{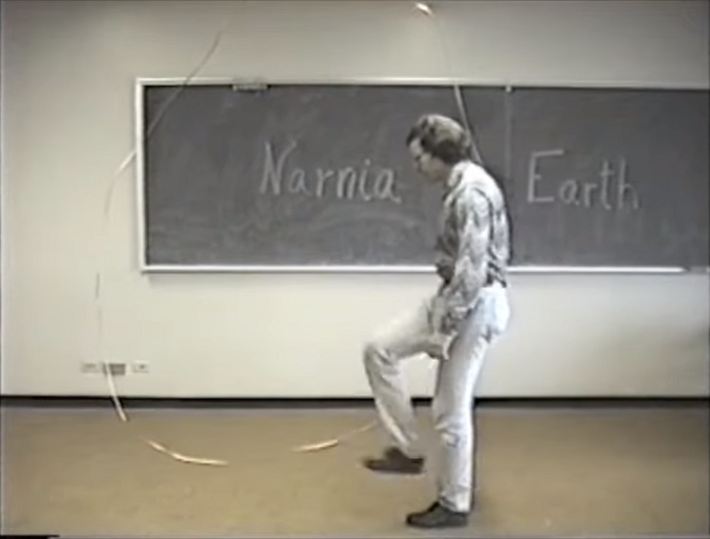}
			\caption{Thurston stepping through a portal generated by the unknot from Earth to Narnia.}\label{fig:thurston}
			\endminipage\hfill
			\minipage[t]{0.48\textwidth}
			\includegraphics[width=\linewidth]{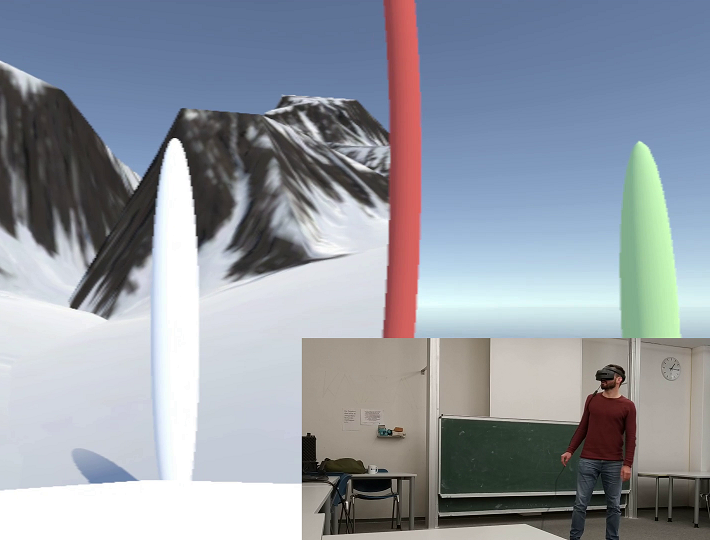}
			\caption{The author stepping through a portal given by the unknot. Screenshot from \url{https://youtu.be/Pgmfsl1e_vA}}\label{fig:kpscreenshot}
			\endminipage\hfill
		\end{figure}

		In a video titled \enquote{Knots to Narnia}~\cite{Thurston2012}, Bill Thurston presents an approach to \enquote{visualize} the cyclic branched cover of a knot by interpreting the knot as a portal to other universes.\footnote{The video was recorded by Tony Phillips as he asked topologists to do \enquote{demos} with knots. To his knowledge, Thurston was the first to illustrate this phenomenon of a branched world in this way.} He demonstrates this using a wire to create different life-sized knotted portals. The wire is \enquote{magical} and, when its ends are joined, creates a \enquote{rip in the fabric of the universe,} creating a portal from our world to a parallel world called \enquote{Narnia} in reminiscence of the novels by C.S. Lewis. The only rule governing the portal is that by circling around the boundary curve twice, one returns to the original world one started in. He then proceeds to explain the phenomena arising in the context of such portals by walking through this wire portal (see Fig.~\ref{fig:thurston}).
		
		This notion of a portal being generated by a ring-shaped object is a quite common theme in movies and videogames, and is mathematically quite simple. Thurston then proceeds to ask a question: What if the wire generating the portal was to be knotted? This leads to different regions in the knot, generating multiple portals. But how many different portals would be generated, and in how many worlds would they lead?
		
		\begin{figure}[H]
			\minipage[t]{0.32\textwidth}
			\includegraphics[width=\linewidth, trim= 50 0 50 0, clip]{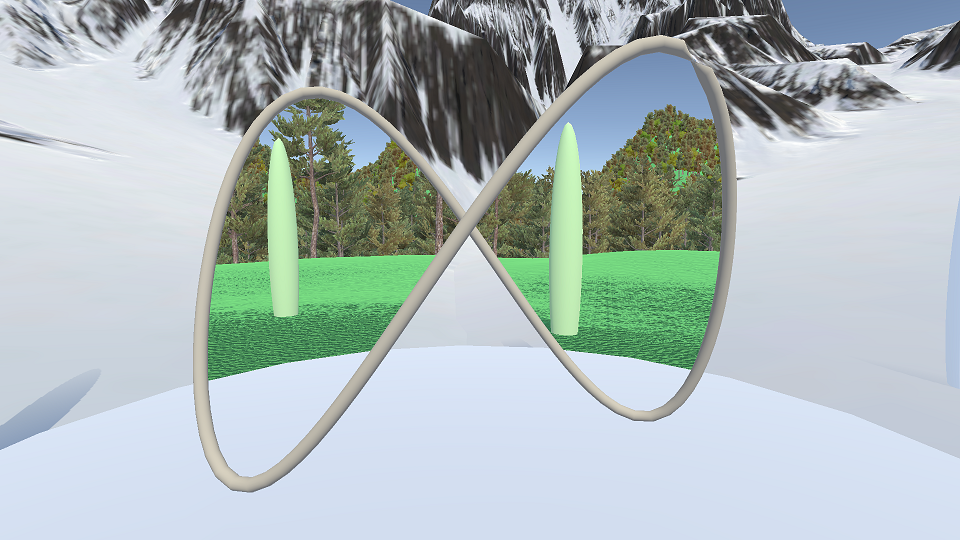}
			\caption{A twisted unknot in \textsc{KnotPortal}, showing two portals into a different world, as seen from the first (ice) world.}\label{fig:twisted1}
			\endminipage\hfill
			\minipage[t]{0.32\textwidth}
			\includegraphics[width=\linewidth, trim= 50 0 50 0, clip]{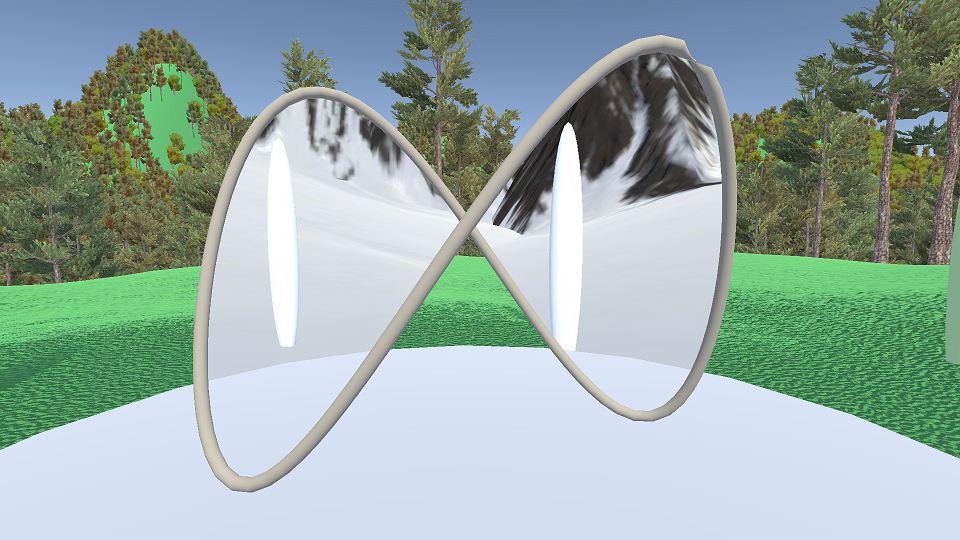}
			\caption{The same twisted unknot, now seen from the second (forest) world.}\label{fig:twisted2}
			\endminipage\hfill			
			\minipage[t]{0.32\textwidth}
			\includegraphics[width=\linewidth, trim= 50 0 50 0, clip]{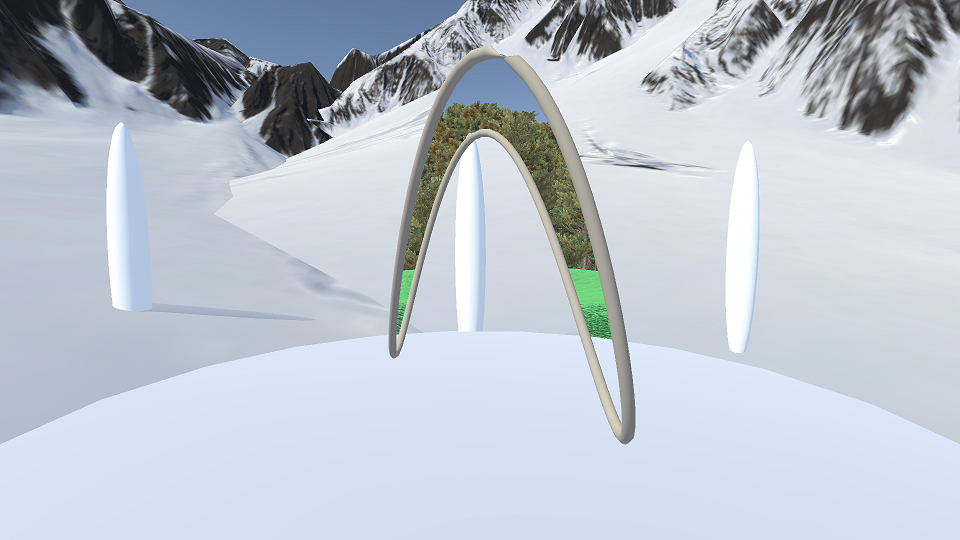}
			\caption{A sideways view of the twisted unknot, revealing why \enquote{both} portals must lead into the same world; there is in fact only one portal.}\label{fig:twistedside}
			\endminipage\hfill
		\end{figure}
		
\begin{figure}[H]
			\begin{subfigure}[t]{0.32\textwidth}
				\centering
				\includegraphics[width=\linewidth]{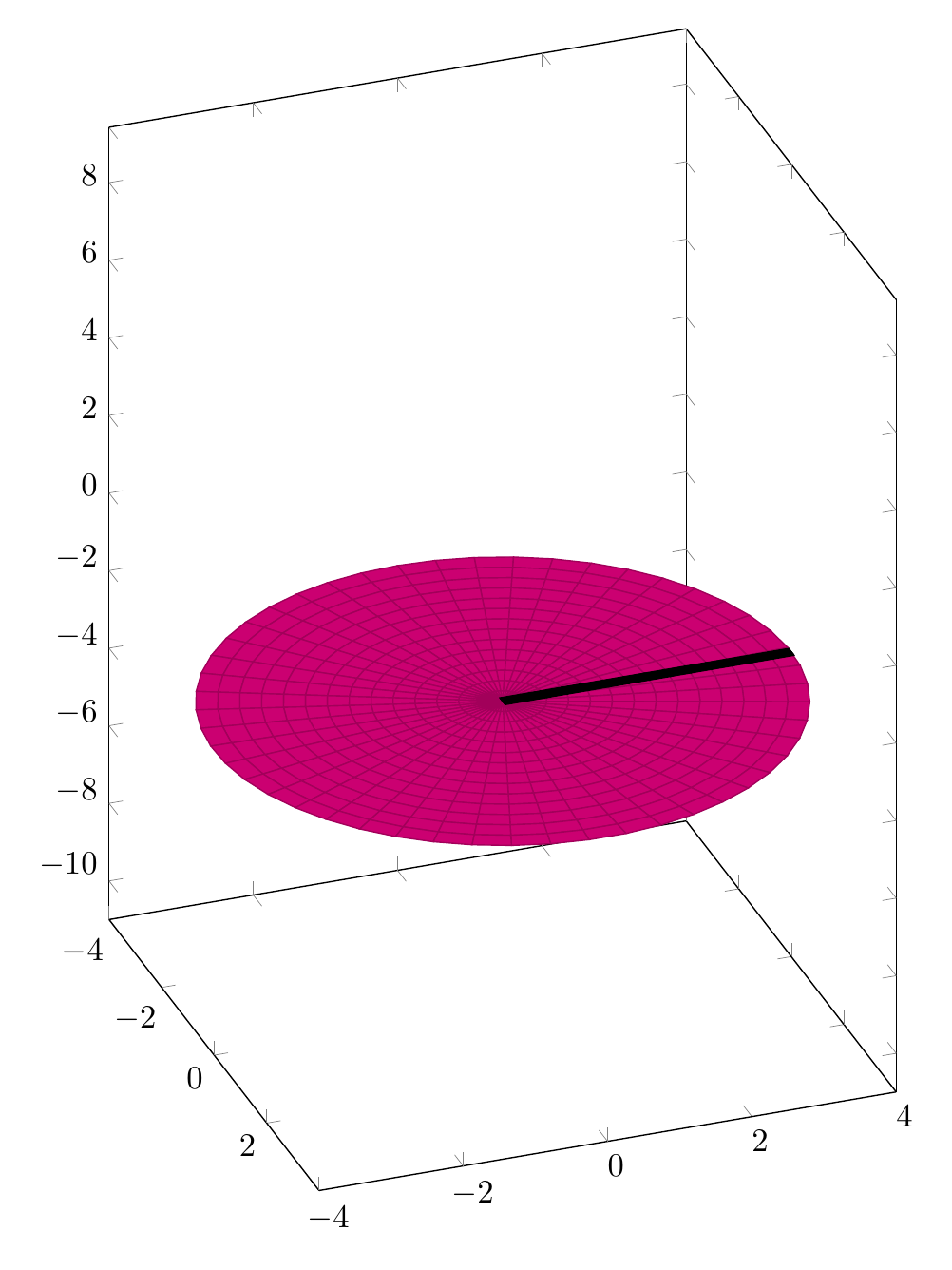}
			\end{subfigure}
			\hfill
			\begin{subfigure}[t]{0.32\textwidth}
				\centering
				\includegraphics[width=\linewidth]{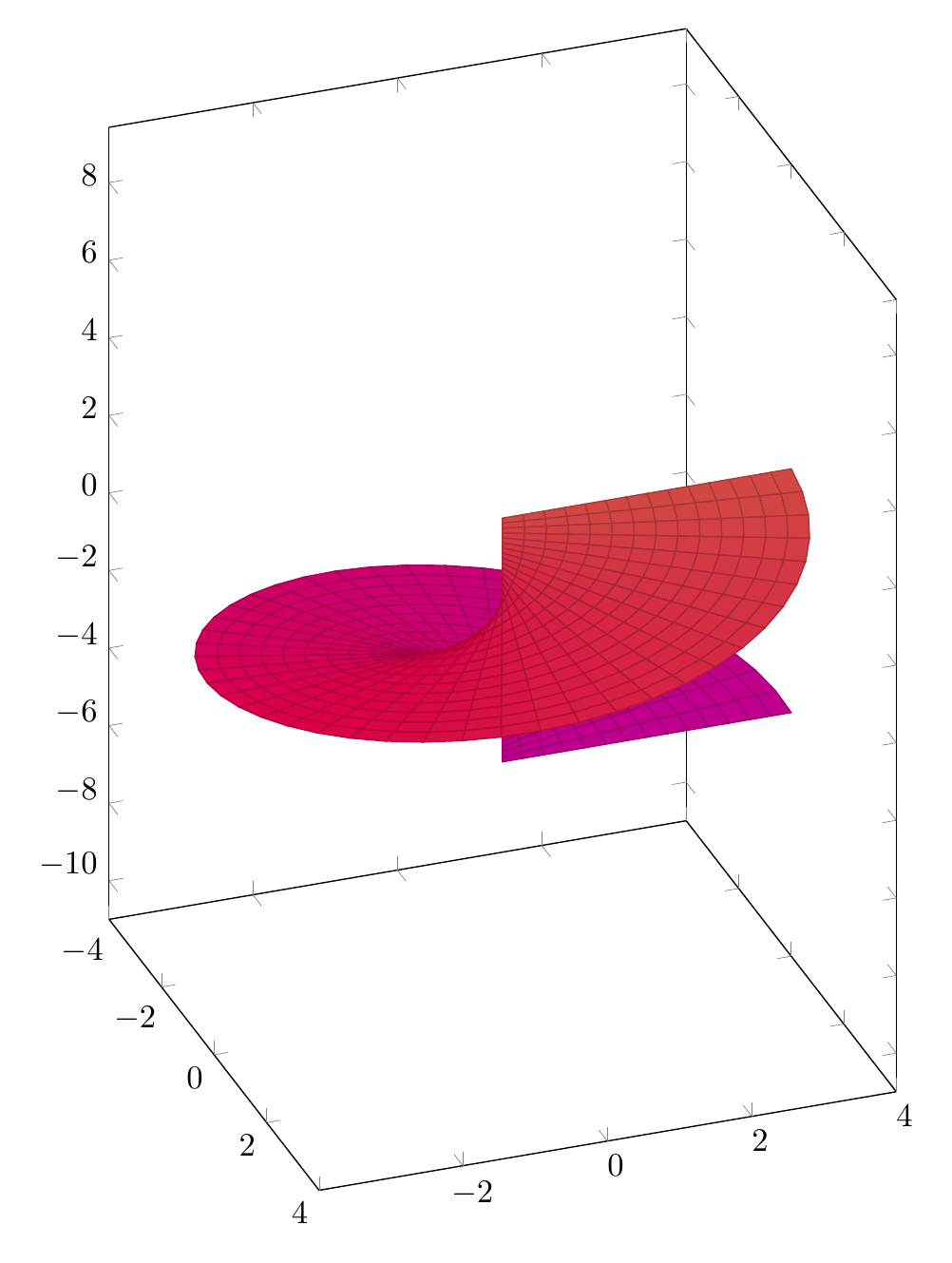}
			\end{subfigure}
			\hfill
			\begin{subfigure}[t]{0.32\textwidth}
				\centering
				\includegraphics[width=\linewidth]{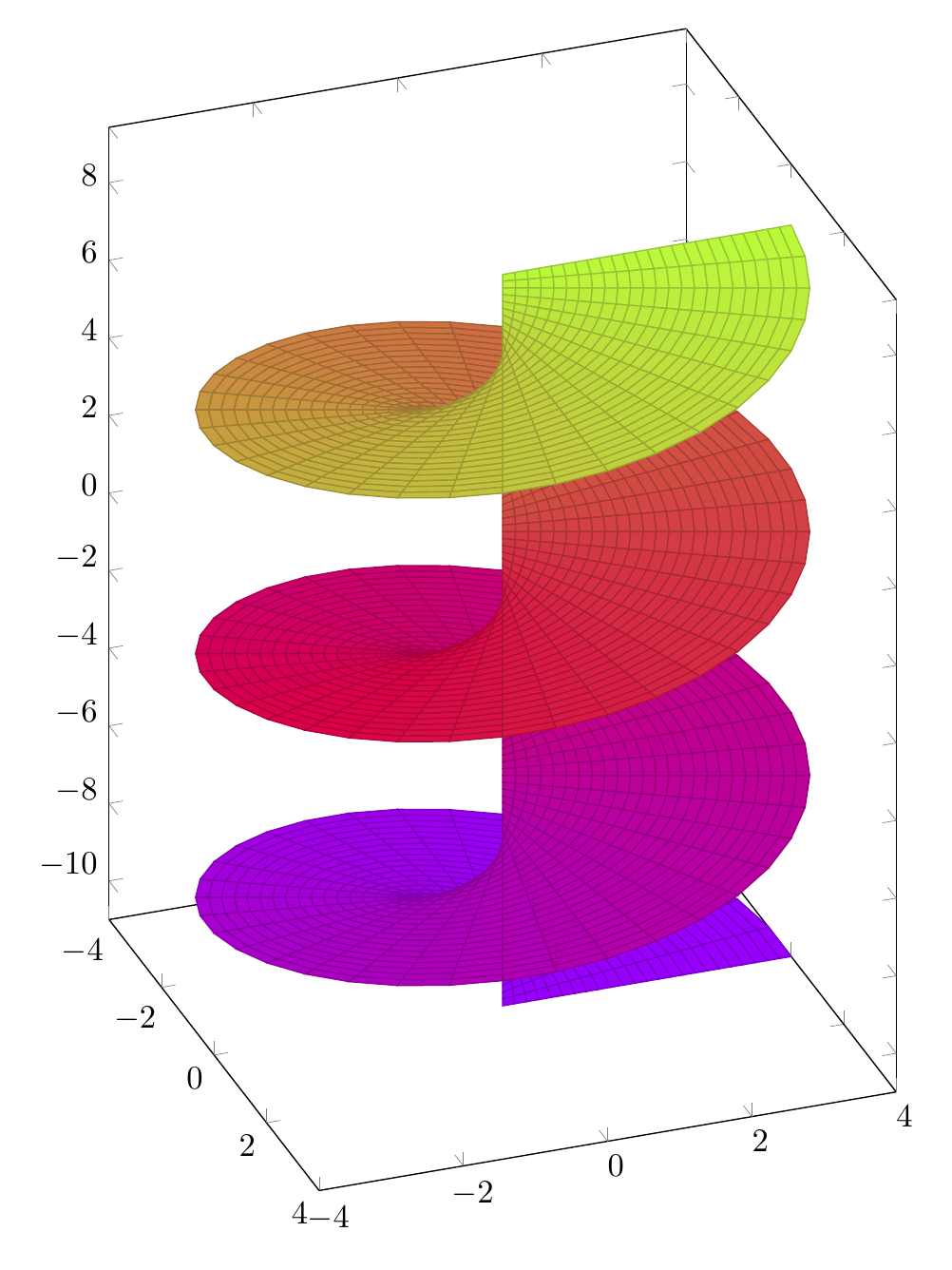}
			\end{subfigure}
			\caption{Creating a branched cover of the complex plane by first cutting from the branch point to infinity, and then gluing together copies of the cut surface along the cutting line.}\label{fig:logCutting}
		\end{figure}

		The object being studied is a cyclic branched cover of order 2. This means that a knot defines a gluing of several sheets of $ \mathbb{R}^3 $, by regarding it as a branching \emph{curve}. Each world is cut along \emph{surfaces} generated by the knot in a way specified in Sec.~\ref{sec:software}, and then glued together according to permutations subject to certain rules. This is analogous to the two-dimensional case, where one has branching \emph{points} and cut \emph{lines} in the construction of, for example, the complex logarithm (see Fig.~\ref{fig:logCutting}).\footnote{This is also an explanation for the common cartoon trope \enquote{behind a stick}, where a character vanishes by running around a tree. It is also what a \enquote{portal} in two-dimensional \enquote{Flatland} would look like.}
		
		This representation of branched covers of knots is fascinating, and for the unknot, it is easy enough to imagine.\footnote{Although it is not completely trivial: If you step through the portal defined by the unknot, and turn around, what do you see?} If, however, the branching curve is knotted, it requires quite a lot of imagination to be able to picture these portals, even for simple cases. This gave the motivation to implement this vision as a computer program, to further recreate Thurston's experience of being able to step through portals as a virtual reality software, giving users the possibility to not only \emph{see} these portals but actually be able to \emph{walk through} them as Thurston did.
		
		In this paper, we describe the implementation of this software and a description of the mathematics involved in the construction of the portals as well as the group structures given by them.
	\section{How to read this article}
		Sec.~\ref{sec:project} gives details of previous work in recreating Thurston's idea. Sec.~\ref{sec:mathematics} contains a short introduction into branched coverings with some interesting examples. In Sec.~\ref{sec:software}, the software \textsc{KnotPortal} is described in detail. Finally, Sec.~\ref{sec:examples} provides examples of branched coverings and the corresponding deck transformation groups.
		
		Readers only interested in the mathematical background of branched coverings of knots need only read Sec.~\ref{sec:mathematics} and maybe \ref{sec:examples} for some examples. For understanding the project, all sections should be read in order, jumping to the examples in Sec.~\ref{sec:examples} on occasion. This last section is of particular interest to those wanting to add own knots to \textsc{KnotPortal}, as it gives an algorithm for doing so.
		
		Regardless the motivation, the reader is strongly advised to try out the software, or at least watch videos of its use, at \url{https://imaginary.org/program/knotportal}.
	\section{Project history}\label{sec:project}
		Previous attempts to model branched covers of knots include the software \enquote{Polycut} by Ken Brakke~\cite{Brakke}. This software was designed \enquote{for visualizing multiple universes connected by a certain kind of wormhole,} with the purpose of illustrating \enquote{the author's contention that soap films are best viewed as minimal cuts in covering spaces.} In the software, the user can view different knots and links and some of their branched covers as differently colored regions, as well as soap films, which are the minimal surfaces separating the sheets.
		
		We wanted to achieve something different, as our goal was to give a real \enquote{world} instead of just colors, as well as to realize a virtual reality experience.\\
		
		There was an attempt to achieve this by porting Ken Brakke's code to CAVE virtual reality technology by George Francis, Alison Ortony, Elizabeth Denne, Stuart Levy and John Sullivan during the illiMath2001 research program, however, this attempt remained unfruitful: \enquote{Though a complete solution to this visualization problem still eludes us, extensive geometrical documentation and evaluation of extant software was undertaken this summer and presented as a PME talk at MathFest, Madison, WI.}, as reported at \url{http://new.math.uiuc.edu/oldnew/im2001/}\\
		
		In this project, we achieved our goal through a new software called \textsc{KnotPortal}, by using the combination of a game engine and a head-mounted virtual reality device capable of room-scale tracking (see Fig.~\ref{fig:kpscreenshot}). In our software, the user can move around in a fully immersive experience featuring different real worlds. It is adaptable as new knots can easily be added, and a non-VR version for use with a normal desktop computer can be used if a VR-headset is not available.

%
	\section{Mathematics background}\label{sec:mathematics}
	\subsection{Branched coverings}
	While this section gives a short overview on branched coverings, interested readers in this topic might want to consult a more comprehensive resource. Most standard textbooks on algebraic topology will do.
		
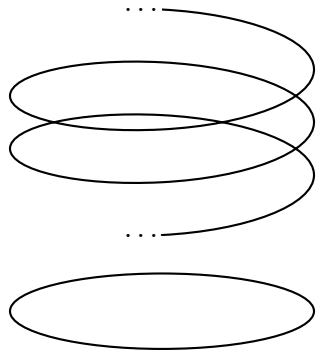
\begin{figure}[H]
		\centering
		\begin{tikzpicture}
		\draw[thick] (0,0) ellipse (2cm and 0.5cm);
		\draw[thick,decoration={aspect=0.31, segment length=7mm, amplitude=2cm,coil},decorate] (0,4) --(0,1);
		\node[text width=1cm] at (0,1) {\ldots};
		\node[text width=1cm] at (0,4) {\ldots};
		\end{tikzpicture}
		\caption{The real line $ \mathbb{R} $ coiled over the sphere $ \mathbb{S}^1 $. Every neighborhood of a point on the sphere has countably infinite many copies above it. As the covering space is simply connected, the cover is called \emph{universal} as it covers every other cover of the base space.}
		\label{fig:expspiral}
	\end{figure}

	A \emph{covering map} is a map $ p $ from a \enquote{covering space} $ E $ to a \enquote{base space} $ X $, such that for any $ x \in X $, the pre-image $ p^{-1}(U_x) $ of any neighborhood $ U_x $ of $ x $ is a disjoint union of open sets $ {\tilde{U}}_{i \in I} $, with $ \tilde{U}_i $ homeomorphic to $ U_x $ for every $ i \in  I $. The cardinality of the index set $ I $ is also called the degree of the cover. In words, this means that every part of the base space has several copies of itself \enquote{above} it. Besides the trivial covering of the disjoint union of copies of a space covering the space itself, the classical example is the \enquote{exponential spiral}. It is defined by the covering map $ p:\mathbb{R} \rightarrow \mathbb{S}^1 $ from the covering space $ \mathbb{R} $ to the base space $ \mathbb{S}^1 $, $ p(t) = \exp^{2\pi i t} $, see Fig.~\ref{fig:expspiral}.\\
	
	If the assumption of \emph{every} point being covered as described above is relaxed to \emph{most} points, one obtains \emph{branched covering maps}. To be precise, a map $ p $ is a branched covering map if it is a covering map for all points but those in a nowhere dense set $ S \subseteq B $, the set of \emph{branch points}. Here, classical examples are the complex square root function as a non-trivial cover of $ \mathbb{C} $ on itself, which is a double, i.e. degree two, covering everywhere but at $ 0 $. Another example is the complex logarithm used as a countably infinite cover of the complex plane, giving rise to the \enquote{logarithmic spiral} in Fig.~\ref{fig:logCutting}. 
	
\begin{figure}[H]
		\centering
		\begin{subfigure}[t]{0.475\textwidth}
			\centering
			\includegraphics[width=\linewidth]{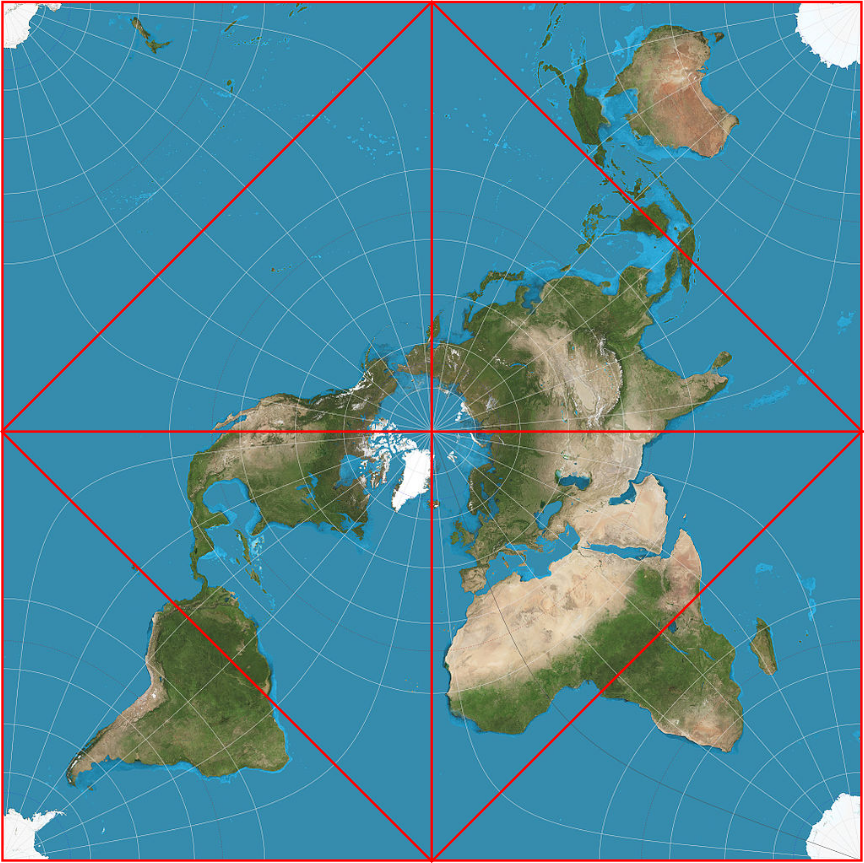}
			\caption{The map projection into a square. The red lines indicate the edges of the octahedron, which is obtained by folding the triangles at the corners backwards. Image from \cite{Strebe2012}.}
			\label{fig:peircesingle}
		\end{subfigure}
		\hfill
		\begin{subfigure}[t]{0.475\textwidth}  
			\centering 
			\includegraphics[width=\linewidth]{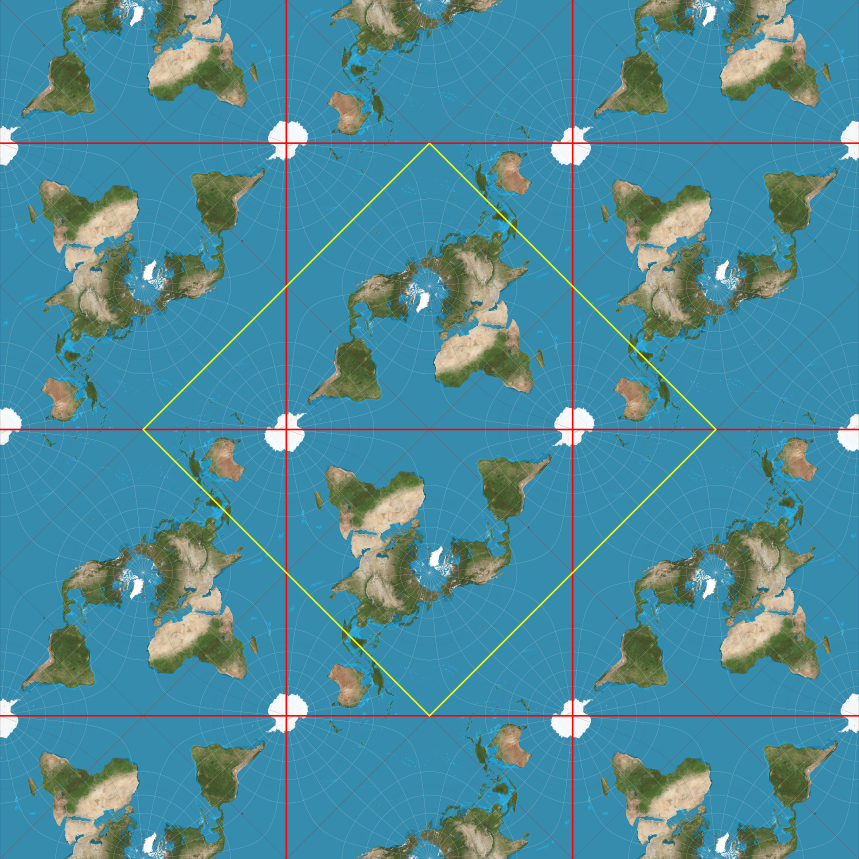}
			\caption{Tiling of the plane by the maps, which are marked in red.}
			\label{fig:peircetile}
		\end{subfigure}
		\vskip\baselineskip
		\begin{subfigure}[t]{0.475\textwidth}   
			\centering 
			\includegraphics[width=\linewidth]{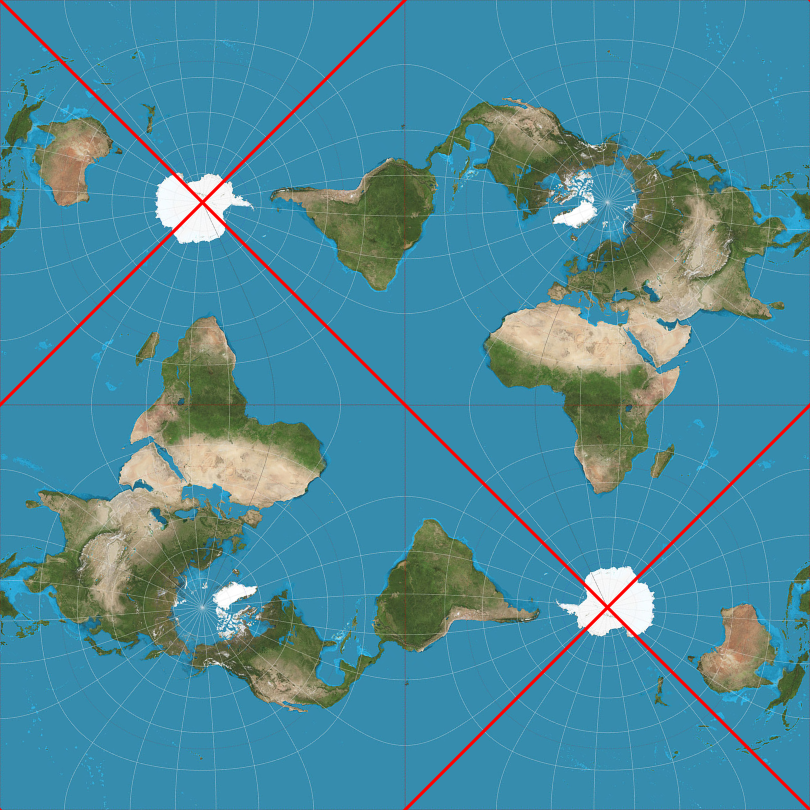}
			\caption{The yellow fundamental region of the torus from Fig.~\ref{fig:peircetile}, with identical corners and identical opposing sides. The four branch points are now at the corner, the center, and the two midpoints of the sides.}
			\label{fig:doubleworld}
		\end{subfigure}
		\quad
		\begin{subfigure}[t]{0.475\textwidth}   
			\centering 
			\includegraphics[width=\linewidth]{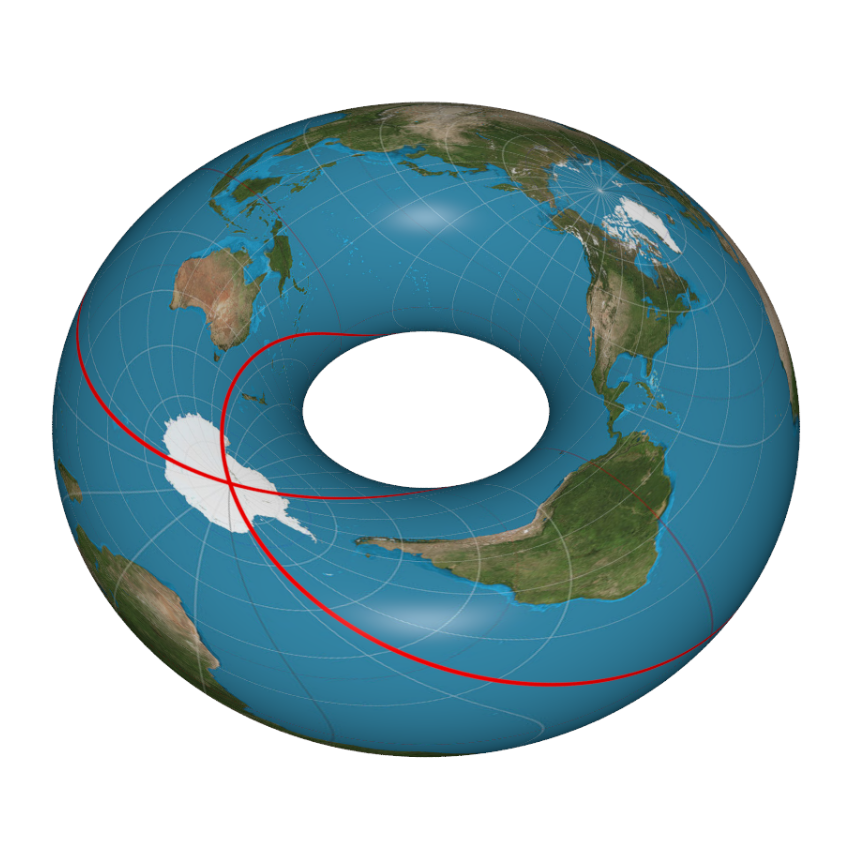}
			\caption{The fundamental region wrapped on a torus. The red lines both connect the Antarcticas. Image generated with \cite{PVR}.}
			\label{fig:peircetorus}
		\end{subfigure}
	\caption{Different views of the Pierce quincuncial map projection.}
	\label{fig:peirce}
	\end{figure}

	Another interesting example is the branched covering of the sphere by the torus, which is also known as Peirce quincuncial projection (see Fig.~\ref{fig:peirce}, or consult \cite{Baez2006} for a complete explanation). It is obtained by projecting the sphere to an octahedron, and then unfolding the octahedron by cutting all edges adjacent to a vertex on the square equator. The flattened version gives a square with the south pole at all corners (See Fig.~\ref{fig:peircesingle}). This square can tile the plane by point reflection on the midpoint of the sides, as shown in Fig.~\ref{fig:peircetile}. This then defines a branched double covering of the sphere by the torus, depicted in Fig.~\ref{fig:peircetorus}, with covering space the torus and base space the sphere.
	
	Every point on the globe is present on the torus twice, except for the branch points, which are only present once. Going around one of the branching points in the covering space also means going around the point on the globe twice. This is not apparent on the map, as Peirce placed the branch points in oceans, making them less visible. 
	
	Branched covers of the sphere are ubiquitous, as made precise by the Riemann existence theorem: every Riemann surface is the branched cover of the sphere~\cite{Harbater2015}.
	
	\subsection{History of the relationship between knots and branched coverings}
	Knots are everywhere in our world, and applications of knot theory range from understanding why headphones get tangled spontaneously~\cite{Raymer2007} to phenomena in quantum physics~\cite{Planat2018}. Although knots are found throughout human history, such as the famous Gordian Knot, their modern mathematical study first began in the 18\textsuperscript{th} century by Vandermonde~\cite{Vandermonde1771} and rised together with topology~\cite{Przytycki2007}. The first applications of known mathematical methods to knots came with Poincaré's \emph{Analysis Situs}~\cite{Poincare1895}. Heegaard used topological methods to compute the 2-fold branch cover of the trefoil knot~\cite{Heegaard1898}, but did not use the result to discriminate the trefoil from the unknot, as this now central problem of knot theory was not of interest to him and was only proved by Tietze in 1908 using the fundamental group~\cite[p.~226]{Stillwell1980}. He used the cover to construct \enquote{Riemann spaces,} analog to the construction of Riemann surfaces in one dimension higher~\cite{Stillwell2012}.
	
	Alexander then proved in \cite{Alexander1920} that \enquote{\emph{Every} closed orientable triangulable n-manifold M is a branched covering of the n-dimensional sphere}, an extension of branched coverings of spheres of the Riemann existence theorem. The theory was even further developed when \cite{Hilden1983} provided a \emph{universal knot}, a knot such that every 3-manifold is a branched cover of the sphere with the knot as a branching set.\footnote{For a more complete history, consult \cite{Artal2017}.}
	
	The knot itself came into the center of attention when Wirtinger extended Heegaard's results and, together with his student Tietze, used the construction to compute a presentation of the fundamental group of the knot complement for every knot~\cite{Epple1999}. The knot group is thus a result of considerations of branched coverings of knots.
	\section{Software}\label{sec:software}	
		The software was created with Unity3D~\cite{UnityTechnologies2017}, the virtual reality gear is HP Mixed Reality\footnote{This is not to be confused with augmented reality; Mixed Reality is just the brand name Microsoft has given its virtual reality technology.}. Scripts are in C\# or, for the shaders, in DirectX 9-style HLSL. The deck transformation groups determining the gluing of the worlds as quotients of the respective knot group, as well as the associated multiplication tables were computed with the help of GAP~\cite{GAP4}.
	\subsection{Input}
		As input, the software is given a knot through some parametrization, as well as a group multiplication table which can be generated with GAP. Examples for knot parametrizations together with group multiplication tables are given in Sec.~\ref{sec:examples}. The software further needs a map defining which \enquote{cone segment} (see below) gets assigned to which group element, the generator-to-cone map.
	\subsection{The setting up of the cut surface}
		At the start of the program, the following steps are carried out.
		\begin{enumerate}
			\item Build all needed worlds
			\item Set up a camera in each world, moving and rotating as the player camera moves and rotates.
			\item Let each camera render to a full-screen sized texture, and assign the textures to the post-processing shader.
		\end{enumerate}
	
		\begin{figure}[H]
			\centering
			\includegraphics[width=0.3\linewidth]{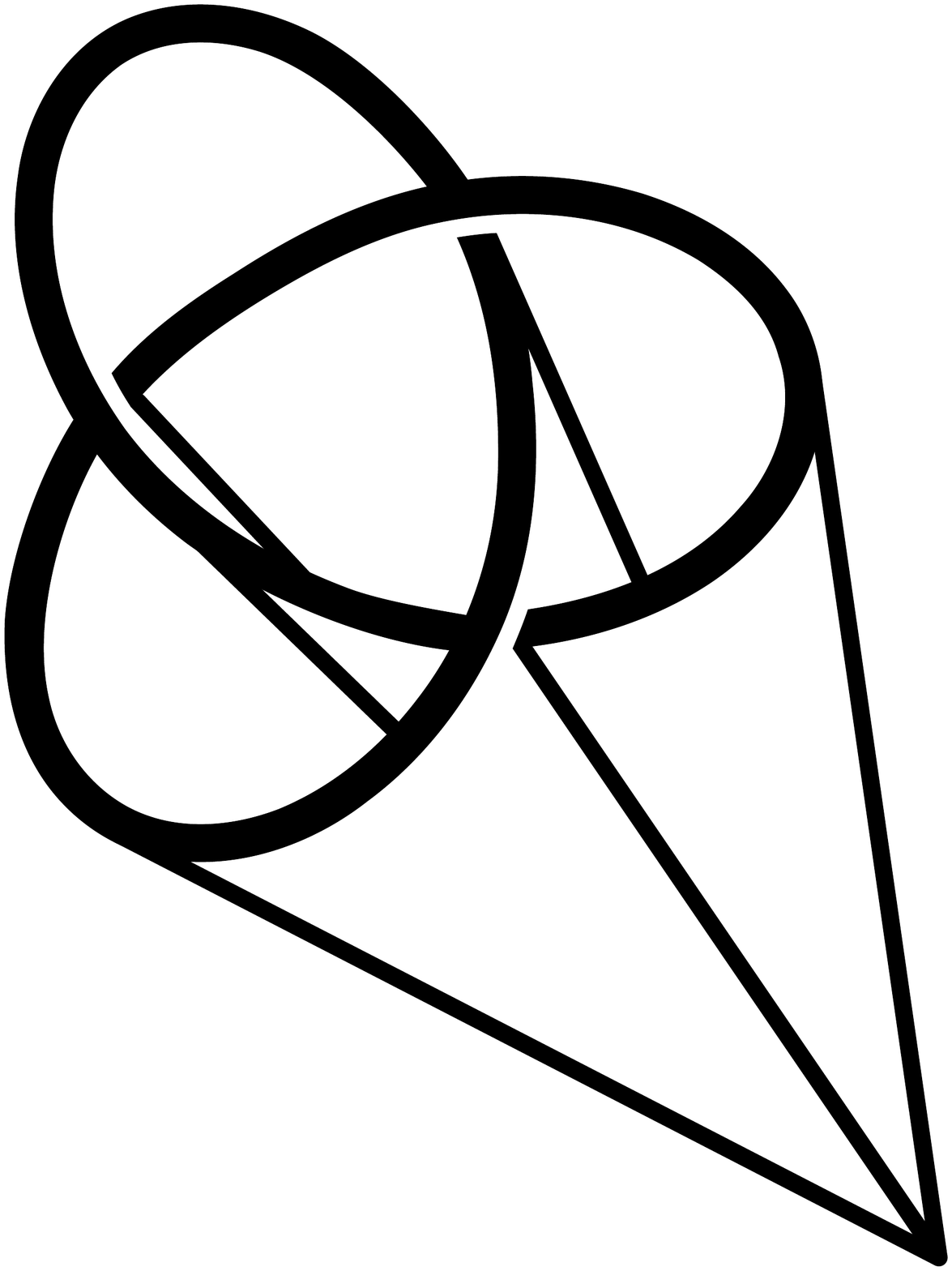}
			\caption[The Heegaard/Reidemeister cone construction for the trefoil knot.]{The Heegaard~\cite{Heegaard1898}/Reidemeister~\cite{Reidemeister1932} cone construction for the trefoil knot. The cone has three self-intersecting lines, resulting in three cone segments (as depicted in \cite{Stillwell1980}).}
			\label{fig:trefoilknotleft}
		\end{figure}
	
		Then, in the first world, we apply the cone construction from \cite{Heegaard1898} to the knot, see Fig.~\ref{fig:trefoilknotleft}. The goal is to provide a cut surface for the gluing of the worlds. This is analogous to the cut \emph{line} given in the construction of the domain of the complex logarithm in Fig.~\ref{fig:logCutting}. In our case, we cut from the branch \emph{curve} to a point at \enquote{infinity} (in the implementation a point sufficiently far away) so that the knot is in general position from its point of view. This defines a cone or cylinder\footnote{also called Reidemeister's cylinder~\cite{Epple1999}} and glue together the different worlds along the cutting \emph{surface}.
		
		\begin{enumerate}
			\item The knot is placed in the world as a tubular mesh around a Catmull-Rom non self intersecting closed spline, given the control points from the discretized parametrization. 
			\item A point $ p $ is chosen, from which a normal knot projection is obtained.
			\item A cone is built from this point by building a mesh formed by the triangles obtained through filling all line segments from $ p $ to every start and end of the line segments of the knot. This results in a sort of cone, possibly self-intersecting.
			\item The cone is cut along the intersections, leading to a number of mesh pieces. These are duplicated and the duplicated has its normals flipped to give a backside.
			\item Each \enquote{cone segment} is assigned a generator of the group according to the provided generator-to-cone map. Its backside gets assigned the inverse of the generator.
		\end{enumerate}
		Now, in each frame, if the knot is visible, perform the following steps on the CPU:
		\begin{enumerate}
			\item Transform the knot's anchor points from world space into screen space.
			\item Using the line segments, divide the screen space into polygonal regions by an algorithm of \cite{Berg2008}.
			\item Find a central point in each region using a C\# port of the \enquote{polylabel} algorithm from \url{https://github.com/mapbox/polylabel} to find the pole of inaccessibility of the region.
			\item Raycast each point from the camera, multiplying the current world generator with every generator from a cone segment encountered along the way. In this way, build a map assigning a generator to each polygonal screen region.
		\end{enumerate}

		\begin{figure}[H]
			\minipage{0.32\textwidth}
			\includegraphics[width=\linewidth]{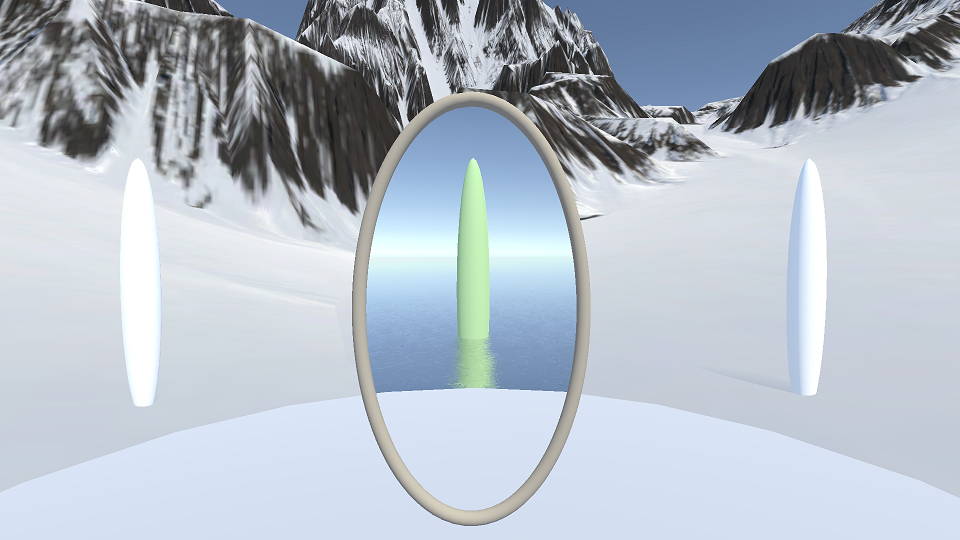}
			\endminipage\hfill
			\minipage{0.32\textwidth}
			\includegraphics[width=\linewidth]{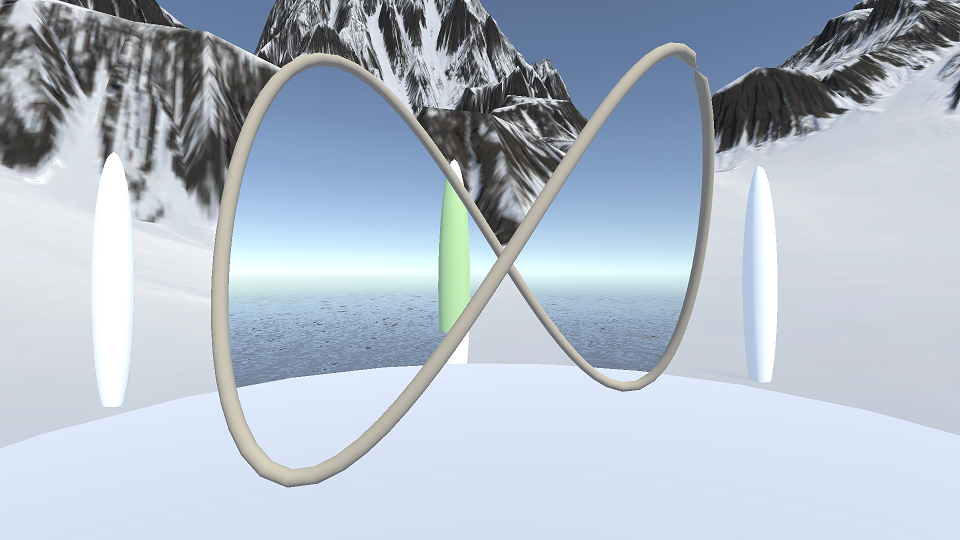}
			\endminipage\hfill
			\minipage{0.32\textwidth}%
			\includegraphics[width=\linewidth]{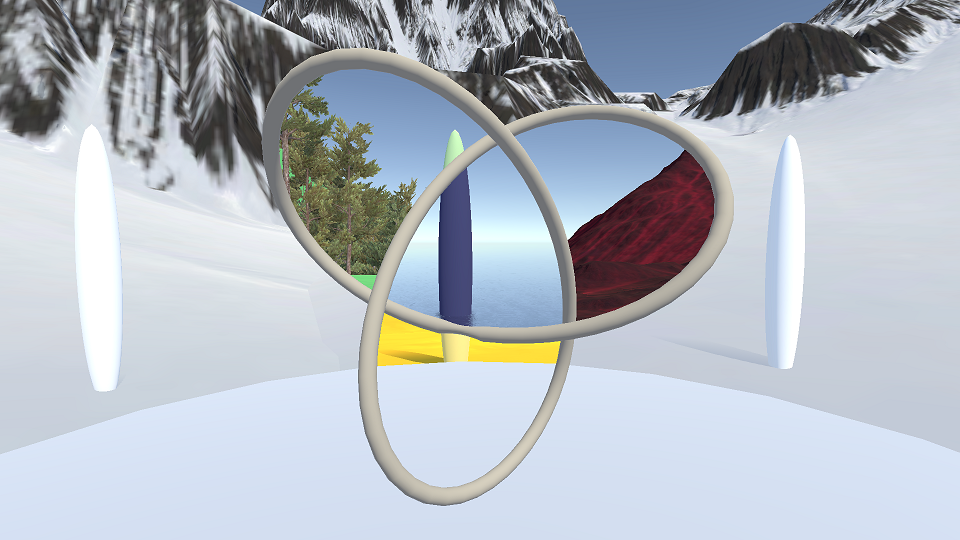}
			\endminipage
			\caption{Portals through the unknot, the twisted unknot, and the trefoil knot}\label{fig:screenshotKnots}
		\end{figure}
		Then run the following steps in the post-processing shader:
		\begin{enumerate}
			\item For each pixel, perform an optimized\footnote{Optimized by first checking if the pixel lies in a bounding box around the polygon, or in a circle of small enough radius around the pole of inaccessibility of the region.} point-in-polygon test.
			\item Assign the pixel the pixel from the camera texture of the world corresponding to the polygon's generator.
		\end{enumerate}
	\subsection{Player teleportation}
		In each frame, perform a raycast from the players old position to his new one. Multiply the current world generator with every cone segment's generator encountered by the raycast, giving the new world. Teleport the player to the point in the same place, but the new world.
		
		This implies that in contrast to expectation, teleportation occurs much later (or earlier, depending on the direction of approach to the knot) as one might think. It does not happen as one \enquote{passes through the portal,} but as one passes through the cut surfaces, i.e. the cone segments, which are the \enquote{real} portal.
	\subsection{World design}
		The software comes with a two different sets of worlds, \emph{simple} and \emph{real} ones. The simple worlds are featureless colored places to enable low-end hardware to run the program, and for a more minimalist experience. 
		
		The other kind are the real worlds (such as in Fig.~\ref{fig:kpscreenshot}), which give the more rich experience. They were designed with several goals in mind. Firstly, they should be interesting enough to give the user a real motivation to step through the portal and look into other worlds. Secondly, they should not be too interesting, as to keep the focus of the experience on the knot and the portals, and not the world. The worlds are also color-coded, to enable the user to speak about \enquote{the white world} or \enquote{the blue world,} which is also helpful in keeping the worlds apart, as well as easing the transition between simple and real worlds. The color codes where taken mainly from naturally occurring colors, with the addition of some colors not present on this planet but possible on other ones~\cite{Kiang2007}.
	
	\section{Example cases}\label{sec:examples}
	These cases all describe branched covers of order 2, i.e. the knot as the branching curve has order 2. So a path going around a knot segment twice is back in the same world (sheet) it started in.
	
	In general, the construction of the deck transformation groups is well-known. Given a (based) cyclic branched covering $ p:(E,e_0) \rightarrow (X,x_0) $ the deck transformation groups can be computed through the Wirtinger presentation together with the fundamental theorem of covering spaces.
	
	The Wirtinger presentation gives the generators of the knot group as loops around the knot strands, together with relations between them for every crossing of the strands.
	
	The fundamental theorem then states that the deck transformation group is isomorphic to
	\begin{equation}
		\pi_1(X,x_0) / p_\star(\pi_1(E,e_0))
	\end{equation}. Given a presentation
	\begin{equation}
		\langle g_1, \ldots, g_n \mid R_1, \ldots, R_m \rangle
	\end{equation} of the knot group, as the covering is cyclic, we have
	\begin{equation}
		p_\star(\pi_1(E,e_0)) \cong \langle g_1^{k_1}, \ldots, g_n^{k_n} \mid R_1, \ldots, R_m \rangle
	\end{equation} for some coefficients $ k_1, \ldots, k_n $. As we restrict ourselves to branched covers of order 2, the coefficients are all 2.
	\subsection{Unknot}
	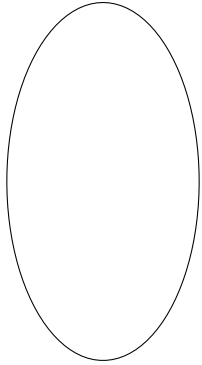
\begin{figure}[H]
		\centering
		\begin{tikzpicture}[scale=1]
		\begin{axis}[
		grid=none,
		axis lines=none,
		axis equal,
		]
		\addplot [domain=0:2*pi,samples=100]({0.8*sin(deg(x))},{1.5*cos(deg(x))}); 
		\end{axis}
		\end{tikzpicture}
		\caption{The unknot under the z-projection}
	\end{figure}
		For the unknot $ K $, the knot group is $ \pi_1(\mathbb{S}^3 \setminus K) $ which is $ \pi_1(\mathbb{S}^1 \times D^2) \cong \mathbb{Z} $ with presentation $ \left\langle a \right\rangle $. Taking the quotient of this group and the subgroup $ \langle a^2 \rangle $, which is the induced by the fundamental group of the covering space, as the simple generating loop has to go around the unknot twice before returning to the basepoint. This results in the presentation $ \left\langle a \mid a^2 \right\rangle $. This is thus a two-fold covering with deck transformation group $ \mathbb{Z}_2 $, or equivalently the (Coxeter) group $ A_1 $.
		
		The unknot is represented in the software through the parametric equations
		$$ \begin{pmatrix}
			0.8 \sin t\\
			1.5 \cos t\\
			0
		\end{pmatrix} $$
		, generates $ |A_1| = 2 $ worlds, and has 1 portal. The group multiplication matrix of $ A_1 $ is
		$ \begin{pmatrix}
			e & a \\ 
			a & e
		\end{pmatrix} $.
		As the cone associated to this knot has no self-intersections, the generator-to-cone map is trivial, assigning every cone segment the group element $ a $.
	\subsection{Twisted Unknot}
	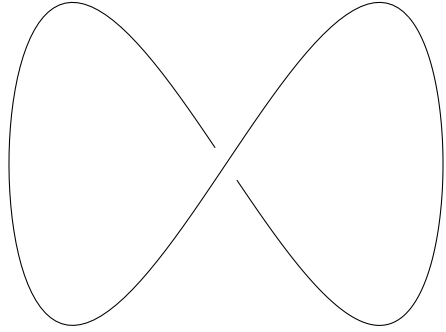
\begin{figure}[H]
		\centering
		\begin{tikzpicture}[scale=1]
		\begin{axis}[
		axis equal,
		grid=none,
		axis lines=none
		]
		\addplot [domain=0:pi-1.05,samples=100]({2*sin(deg(x+1))},{3*cos(deg(x+1))*sin(deg(x+1))});
		\addplot [domain=pi-0.95:2*pi,samples=100]({2*sin(deg(x+1))},{3*cos(deg(x+1))*sin(deg(x+1))});  
		\end{axis}
		\end{tikzpicture}
		\caption{The twisted unknot under the z-projection}
	\end{figure}
		This case is of course the same as the unknot from a knot theoretical standpoint.
		
		As for the implementation, the knot is given by
		\[ \begin{pmatrix}
		2 \sin (t+1)\\
		3 \sin (t+1)  \cos (t+1)\\
		\sin t
		\end{pmatrix} \], but as there are two portals leading to the same world, the generator-to-cone map assigns $ a $ to both cone segments.
	\subsection{Trefoil knot}
	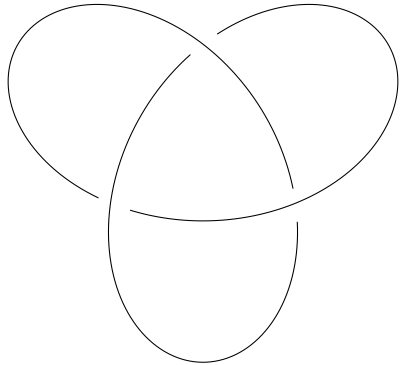
\begin{figure}[H]
		\centering
		\begin{tikzpicture}[scale=1]
		\begin{axis}[
		axis equal,
		grid=none,
		axis lines=none
		]
		\addplot [domain=0:pi-1.37,samples=100]({sin(deg(x))+2*sin(deg(2*x))},{cos(deg(x))-2*cos(deg(2*x))});
		\addplot [domain=pi-1.27:pi+0.72,samples=100]({sin(deg(x))+2*sin(deg(2*x))},{cos(deg(x))-2*cos(deg(2*x))});
		\addplot [domain=pi+0.82:2*pi-0.31,samples=100]({sin(deg(x))+2*sin(deg(2*x))},{cos(deg(x))-2*cos(deg(2*x))});
		\addplot [domain=2*pi-0.21:2*pi,samples=100]({sin(deg(x))+2*sin(deg(2*x))},{cos(deg(x))-2*cos(deg(2*x))});
		\end{axis}
		\end{tikzpicture}
		\caption{The trefoil knot under the z-projection}
	\end{figure}
		For the trefoil knot $ K $, the knot group is $ \left\langle a,b \mid a^3 = b^2 \right\rangle $ as the trefoil knot is the $ (2,3) $ torus knot~\cite{Stillwell1980}. Alternatively, it can be given by $ \left\langle x,y \mid xyx = yxy \right\rangle $~\cite[p.~61]{Rolfsen2003}. By using $ xyx = yxy \cong xyxxyx = xyxyxy \cong yxyxyx = (xy)^3 \cong yxxyxx = (xy)^3 \cong (yxx)^2 = (xy)^3$ we can see the isomorphism between the two presentations. Adding the relations $ x^2 $ and $ y^2 $, we obtain the presentation $ \left\langle a,b \mid (xy)^3, x^2, y^2 \right\rangle $. This is the dihedral group of the triangle, and a Coxeter group with Coxeter matrix $ \begin{pmatrix}
		1 & 3  \\ 
		3 & 1
		\end{pmatrix} $. The group order 6 implies the construction of 6 worlds from this knot. In general, the $ r $-fold branched covering of the torus knots of type $ (p,q) $ is a Brieskorn manifold $ M(p,q,r) $, the intersection of the 5-sphere $ \mathbb{S}^5 $ in $ \mathbb{C}^3 $ with the equation given through $ z_1^p + z_2^q + z_3^r = 1$.~\cite{Planat2018}.
		
		In \textsc{KnotPortal}, the trefoil knot is represented through the parametric equations
		\[ \begin{pmatrix}
			\sin t + 2 * \sin 2t\\
			\cos t - 2 * \cos 2t\\
			-\sin 3t
		\end{pmatrix} \]. The group multiplication matrix of $ D_3 $ is 
		\[ \begin{pmatrix}
			e & a & b & c & d & f \\ 
			a & e & d & f & b & c \\ 
			b & f & e & d & c & a \\ 
			c & d & f & e & a & b \\ 
			d & c & a & b & f & e \\ 
			f & b & c & e & a & d
		\end{pmatrix}  \]. The generator-to-cone map assigns the elements $ a $, $ b $, and $ c $ to the three cone segments, respectively.
		
\begin{figure}[H]
			\centering
			\includegraphics[width=0.5\linewidth]{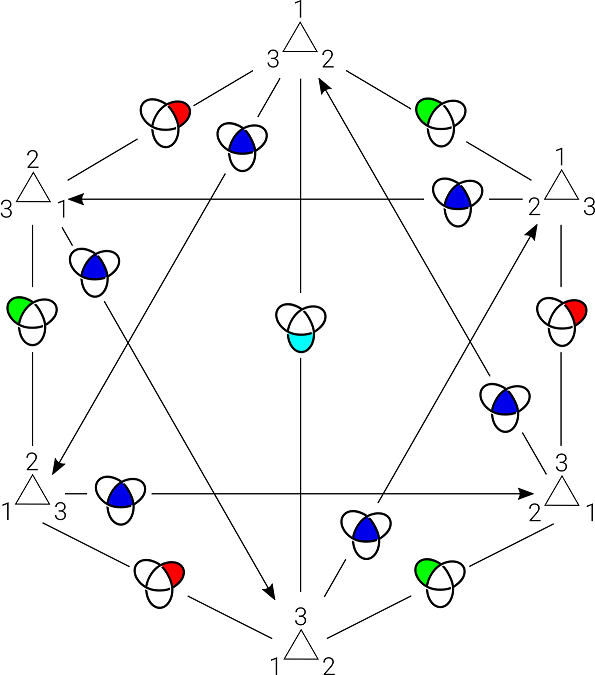}
			\caption{The relationship between the elements of the symmetry group of the triangle $ D_3 $ and the portals generated by the trefoil knot, after the drawing in \cite{Thurston2012}. The outer portals correspond to reflections, the inner portal to a rotation.}
			\label{fig:dihedral}
		\end{figure}

		The relationship between the group and the portals of the trefoil knot is detailed in Fig.~\ref{fig:dihedral}
		
	\subsection{Figure eight knot}
	\begin{figure}[H]
		\centering
		\begin{tikzpicture}[scale=1]
		\begin{axis}[
		axis equal,
		grid=none,
		axis lines=none
		]
		\addplot [domain=0:1.02,samples=100]({(2 + cos(deg(2*x)))* cos(deg(3*x))},{(2 + cos(deg(2*x)))* sin(deg(3*x))});
		\addplot [domain=1.08:2.59,samples=100]({(2 + cos(deg(2*x)))* cos(deg(3*x))},{(2 + cos(deg(2*x)))* sin(deg(3*x))});
		\addplot [domain=2.65:3.635,samples=100]({(2 + cos(deg(2*x)))* cos(deg(3*x))},{(2 + cos(deg(2*x)))* sin(deg(3*x))});
		\addplot [domain=3.695:4.16,samples=100]({(2 + cos(deg(2*x)))* cos(deg(3*x))},{(2 + cos(deg(2*x)))* sin(deg(3*x))});
		\addplot [domain=4.22:2*pi,samples=100]({(2 + cos(deg(2*x)))* cos(deg(3*x))},{(2 + cos(deg(2*x)))* sin(deg(3*x))});
		\end{axis}
		\end{tikzpicture}
		\caption{The figure eight knot under the z-projection}
	\end{figure}
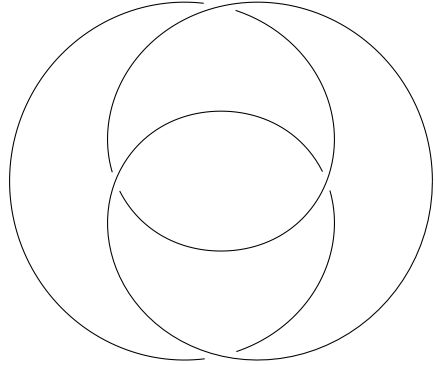
	The presentation of the figure eight knot is $ \left\langle x,y \mid x^{-1}yxy^{-1} = yx^{-1}yx \right\rangle $~\cite[p.~58]{Rolfsen2003}. Again adding the relations $ x^2 $ and $ y^2 $, one obtains $ \left\langle x,y \mid (xy)^5, x^2, y^2 \right\rangle $, which is again a Coxeter group, namely $ H_2 $, which is of order 10. This knot thus generates 10 worlds.
	
	In the software, it is represented through
	\[ \begin{pmatrix}
	(2 + \cos 2t) \cos 3t\\
	(2 + \cos 2t) \sin 3t\\
	\sin 4t
	\end{pmatrix} \]
	The group multiplication table is 
	\[ \begin{pmatrix}
	a & b & c & d & e & f & g & h & i & j\\
	b & a & d & c & f & e & h & g & j & i\\
	c & j & e & b & g & d & i & f & a & h\\
	d & i & f & a & h & c & j & e & b & g\\
	e & h & g & j & i & b & a & d & c & f\\
	f & g & h & i & j & a & b & c & d & e\\
	g & f & i & h & a & j & c & b & e & d\\
	h & e & j & g & b & i & d & a & f & c\\
	i & d & a & f & c & h & e & j & g & b\\
	j & c & b & e & d & g & f & i & h & a
	\end{pmatrix} \]	
	\subsection{Solomon's Seal knot}
	This is the $ (5,2) $-torus knot. Its parametric equation is thus given by~\cite{VonSeggern2016}:
	\[ \begin{pmatrix}
	(3 + \cos 5t) \cos 2t\\
	(3 + \cos 5t) \sin 2t\\
	\sin 5t
	\end{pmatrix} \]	
	and the presentation of its group is $ \left\langle x,y \mid xyxyxy^{-1}x^{-1}y^{-1}x^{-1}y^{-1}\right\rangle $~\cite{Livingston1993}. After adding the relations for the generators, the order two covering group of this knot is thus the same as for the Figure eight knot.
	\subsection{Hopf Link}
	\begin{figure}[H]
		\centering
		\begin{tikzpicture}
		\begin{axis}[
		axis equal,
		grid=none,
		axis lines=none
		]
		\addplot [domain=0:pi-0.58,samples=100]({1*sin(deg(x))},{1*cos(deg(x))});
		\addplot [domain=pi-0.48:2*pi,samples=100]({1*sin(deg(x))},{1*cos(deg(x))}); 
		\addplot [domain=0:pi+2.56,samples=100]({1*sin(deg(x))+1},{1*cos(deg(x))}); 
		\addplot [domain=pi+2.66:2*pi,samples=100]({1*sin(deg(x))+1},{1*cos(deg(x))}); 
		\end{axis}
		\end{tikzpicture}
		\caption{The Hopf link under the z-projection}
	\end{figure}
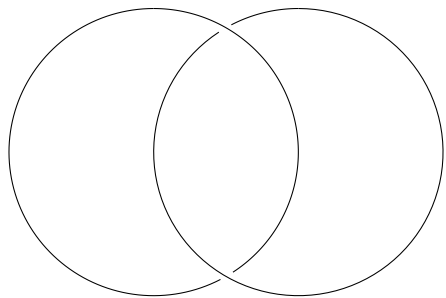
	Each of the branching curves gives a generator, and the two commute, so the deck transformation group is $ \left\langle a,b \mid a^2, b^2, (ab)^2 \right\rangle $. This group is a Coxeter group with matrix $ \begin{pmatrix}
	1& 2\\ 
	2& 1
	\end{pmatrix} $, which is $ \mathbb{Z}_2^2 $, or equivalently, $ A_1^2 $. This results in 4 worlds and 3 portals.\footnote{At the present time, the support of links is not implemented in the software, but could certainly be achieved without much change to the methods.}
		
	\section*{Acknowledgments}
	Thanks to Marc Sauerwein, John Sullivan and Roice Nelson for their valuable advice.
	\bibliographystyle{alpha}
	\bibliography{literatur}
\end{document}